# Assessing the Dissipative Capacity of Particle Impact Dampers Based on their Nonlinear Bandwidth Characteristics

Xiang Li[1,2*], Alireza Mojahed[4], Li-Qun Chen[1*], Lawrence A. Bergman[3],
Alexander F. Vakakis[4]

1. Shanghai Institute of Applied Mathematics and Mechanics, School of Mechanics and Engineering Science, Shanghai University, Shanghai 200444, China
2. Faculty of Civil Engineering and Mechanics, Jiangsu University, Zhenjiang 212013, China
3. Aerospace Engineering, University of Illinois, Urbana, IL 61801, USA
4. Mechanical Science and Engineering, University of Illinois, Urbana, IL 61801, USA

**Abstract**

The dissipative capacity as quantified by the nonlinear bandwidth measure of impulsively loaded linear primary resonators or primary structures (PSs) coupled to particle impact dampers (PIDs) is assessed. The considered PIDs are designed by initially placing different numbers of spherical, linearly viscoelastic granules at different 2D initial topologies and clearances. The strongly nonlinear and highly discontinuous dynamics of the PIDs are simulated via the discrete element method taking Hertzian interactions, slipping friction and granular rotations into account. The general definition of nonlinear bandwidth is used to evaluate the energy dissipation capacity of the integrated PS-PID systems. Moreover, the effect of the dynamics of the PIDs on the time-bandwidth product of these systems is studied, as a measure of their capacity to store or dissipate vibration energy. It is found that the initial topologies of the granules in the PID drastically affect the time-bandwidth product, which, depending on shock intensity, may break the classical limit of unity which holds for linear time-invariant dissipative resonators. The optimal PS-PID systems composed of multiple granules produce large nonlinear bandwidths, indicating strong dissipative capacity of broadband input energy by the PIDs. Additionally, in the optimal configurations, the time-bandwidth product, i.e., the measure of the frequency bandwidth of the input shock that is stored in the PS-PID system, in tandem with the amount of time it takes for the system to dissipate $(1/e)$ of the initial energy, can be tuned either above or below unity by varying the applied shock intensity. The implications of these findings on the dissipative capacity of the system considered are discussed, showing that it can be predictively assessed so that PIDs can act as highly effective nonlinear energy sinks capable of rapid and efficient suppression of vibration induced by shocks.

[*]corresponding authors: xianglichn@outlook.com and lqchen@shu.edu.cn

# 1. Introduction

Topics related to particle dampers or particle impact dampers (PIDs) have attracted much attention [1,2]. In a typical configuration, a primary structure (PS) in the form of a resonator is attached to a rigid container or cavity containing several colliding spherical particles, or granules, designated as the PID. In most cases the granules are assumed to be composed of linear viscoelastic material, whereas the granule-granule and granule-container interactions are modeled through a combination of Hertzian and frictional forces (as detailed in the discussion that follows). It was demonstrated that PIDs can be a useful vibration/shock suppression passive tool, since they are capable of transferring input energy into the granules, and locally and rapidly dissipating it through inelastic collisions and frictional interactions between the granules themselves and the boundary granules and their container walls [3]. Key advantages of PIDs include their relatively simple and modular composition, their insensitivity to temperature fluctuations, and, as shown in [3] and in this work, their broadband vibration/shock suppression capacity.

Modeling the dissipative capacity of PIDs is a key consideration for effective predictive design to achieve desired performance. Recently, Masmoudi et al. [4] used a loss factor approach to predict energy dissipation in PIDs, reporting that large mass or high excitation magnitude leads to more energy dissipation. Lu et al. [5] proposed a nondimensional energy dissipation factor to select material properties, such as modulus of elasticity and yield strength, for better PID energy dissipation. Xiao et al. [6] applied particle dampers in a gear transmission; according to the results of simulations and experiments, they concluded that a small restitution coefficient in granular interactions contributes to better energy dissipation at high rotational gear speeds. Yan et al. [7] found that at low excitation amplitudes and frequencies it is necessary to account for frictional effects between the granules of the PID and its container walls. Niklas and Robert [8] proposed a coupling method based on a combination of a reduced loss factor and an effective particle mass to estimate the overall damping of the structure with a PID attached; they showed that the PID can effectively suppress even multi-modal responses of the PS-PID configuration. According to an equivalent single particle model, Lu et al. [9] computed the root-mean-square responses of a PS coupled to a PID under a broadband random excitation;

they validated that the PID can suppress the broadband random responses. Sack et al. [10] proved the sharp transition between the gas-like and collective motions of the granular dampers in a microgravity environment, and concluded that collective motions favor energy dissipation. Ye et al. [11] developed a deep transfer learning technique with less simulation results and more experimental data to investigate the efficiency of a PID subject to harmonic excitation of the PS over a wide frequency range.

Moreover, alternative novel forms of PIDs were proposed to improve vibration attenuation performance. For example, Zhang et al. [12] designed an electromagnetic particle damper, and confirmed that introduction of the magnetic forces can enhance the vibration dissipative capacity of the PID. Varela-Rosales et al. [13] studied PIDs composed of granules of the same material but with different sizes; numerical simulation results indicated that this "mixture" PID dissipates energy much more efficiently compared to a same-mass traditional PID composed of identical granules. In another work, Gnanasambandham et al. [14] incorporated rigid obstacle grids into the PIDs in order to promote relative motions between granules and, accordingly, obtain efficient energy dissipation at low external excitation levels. Hu et al. [15] extended the study of PIDs with obstacle grids, and reported that these grids effectively improve energy dissipation for a range of the excitation intensities or filling ratios; in addition, they demonstrated directed energy transfer from the PS to PID. Żurawski et al. [16] proposed a kind of adaptive tuned PID with a singular granule and a changeable height of the container, and showed the efficient performance of vibration attenuation by modifying the container height in real-time.

Owing to their strongly nonlinear and highly discontinuous dynamics [2], PIDs are capable of absorbing and locally dissipating vibrational energy from the host PS [1,3,17]. Bai et al. [17] studied piston-type and box-type particle dampers from the view of energy analysis via the discrete element method; they reported that both energy transfer and energy consumption capacity govern the performance of the PIDs. Lu et al. [1] suggested that the damping mechanism of the PIDs can be explained further from the perspective of the nonlinear energy sink (NES), and targeted energy transfer (TET). TET is denoted by the passive directional transfer from a directly excited primary structure to a strongly nonlinear attachment, designated as the NES, via transient or sustained nonlinear resonant captures [18]. Xiang et al. [3] studied

PIDs with topologies to irreversibly transfer shock energy from the PS to the attached PID serving essentially as a granular NES; they also discussed the relationships between TET and granular dynamics.

Yang et al. [19] proposed a NES with the time-varying mass to trigger TET at much lower energy levels and accordingly overcome the well-known critical energy threshold required for "activation" of the NES. Ding et al. [20–22] introduced NESs with piecewise characteristics, showing that this is an efficient way to suppress multi-mode resonances in the PS. Wang et al. [23] integrated a lever-based NES with a magnetostrictive energy harvester, and experimentally confirmed dual functionality of the NES as vibration suppressor and energy harvester. Zeng and Ding [24] designed a tri-stable NES via the pre-compressed beam with magnets, and studied its capacity for vibration suppression of small-magnitude excitations. Liu and Wang [25] added geometrically nonlinear damping in a multi-degree-of-freedom (multi-DOF) NES for the purpose of suppressing large-amplitude responses. Dou et al. [26] applied particle damping into a bistable NES to reinforce its dissipative capability for torsional vibration suppression of rotor systems.

In an alternative approach, Vakakis [27] proposed irreversible low-to-high frequency energy transfer or modal energy redistribution as a robust and efficient vibration control method in systems with or without NES. It is reasonable that low-to-high frequency energy scattering can be an effective dissipative method due to the well-established inherent property of higher structural modes for low-amplitude vibration and enhanced energy dissipation. Additional investigations were carried out to assess the efficacy of this concept [28–33]. Fang et al. [28] introduced multiple vibro-impact NESs to a harmonically-excited cantilever beam for effective low-to-high frequency energy transfer. Li et al. [29] designed and optimized a two-degree-of-freedom (two-DOF) NES with geometrically nonlinear damping; they used a normalized effective damping measure to validate the high efficiency of shock energy transfer from low to much higher modes of a high-rise building. Theurich et al. [30] developed a useful semi-analytical method to investigate a flexible structure coupled to an impact absorber (behaving as a vibro-impact NES), and confirmed that the NES not only absorbs and dissipates energy, but also scatters energy among the modes of the flexible structure. The so-called intermodal targeted energy transfer (IMTET) mechanism [31] related to low-to-high frequency energy scattering within the modal space of a linear PS was successfully applied to rapidly diminish the shock response of a cantilever beam [32], and the seismic response of a tall building structure [33]. These works indicate that passive vibration suppression based on the TET mechanism in the time or frequency domains is an effective approach to quickly attenuate vibration of structures under external excitations. However, *how to exactly measure or quantify the dissipative capability of a PS coupled to an NES or a PID in the time or frequency domains is a separate interesting issue that deserves special consideration*. Indeed, answering this question is vital in order to optimize PIDs for vibration reduction.

The dissipative capacity of a PID may be associated with the *nonlinear bandwidth* of the integrated PS-PID system [34]; vibration analysis in this context emerges as a new and promising research field. Mojahed et al. [35] provided a general definition of nonlinear bandwidth as a way of extending and generalizing the classical notion of bandwidth defined for linear resonators with relatively weak damping. This new definition, referred to also as root-mean-square (RMS) bandwidth, applies to nonlinear/linear and time-invariant/time-variant resonators and overcomes the limitation of the classical half-power bandwidth which is valid only for linear dissipative systems. By proposing a definition based on the root mean square (RMS) bandwidth and the envelope of the decaying energy of a resonator, one recovers the original purpose of bandwidth, which is to quantify the overall dissipative capacity of the system or, equivalently, to describe how localized the energy of the system is in the time and frequency domains. Since the RMS bandwidth of the energy signal is linearly proportional to the inverse of the variance of the energy signal in the time domain (according to the Fourier uncertainty principle), the new bandwidth definition provides an accurate measure of the dissipation rate (capacity) of the free decay of the response of the resonator and is its inherent property [34,35]. Consequently, the concept of nonlinear bandwidth can be used to quantify the dissipative capacity of a PS attached to a PID, overcoming the strongly nonlinear and highly discontinuous dynamics of this system, and the strong tunability of the dynamics with energy.

Additionally, nonlinear bandwidth enables investigation of the *time-bandwidth product* of a nonlinear time-invariant system [36]. Considering a linear time invariant (LTI) resonator, its vibration will decay with time owing to some loss mechanism(s) with a (total) decay rate $\lambda$. Assuming the system is underdamped and at resonance, it holds that $\lambda = \Delta\omega$, where $\Delta\omega$ is the (classical half-amplitude) bandwidth. Introducing the "storage time" $\Delta t = 1/\lambda$ (i.e., the time that a resonator can store the vibration energy), we deduce that the product of the storage time and the bandwidth should be always equal to unity, i.e., $\Delta t\, \Delta \omega = 1 \Rightarrow \Delta t = 1/\Delta \omega$, and this holds for any single-DOF dissipative LTI resonator; this is referred to as the *classical "time-bandwidth (T-B) limit"* [36] for LTI resonant devices. This limit ascertains that an LTI resonator can either exhibit weak dissipative capacity and narrowband resonance, or strong dissipative capacity and broadband resonance, *but not simultaneously both*. Breaking the T-B product, e.g., above unity would imply that we could harness the "best of both worlds" in a resonator system; namely, we could acquire simultaneous capacity for storing vibrations for longer times (i.e., possessing weaker dissipative capacity) *and* at broader bandwidths compared to the linear case. Conversely, tuning this limit below unity would mean that, e.g., for a fixed large bandwidth, the release of a vibration would be faster than what is normally allowed (for that bandwidth). Mojahed et al. [36] theoretically and experimentally proved that the time-bandwidth product of a nonlinear time-invariant system can be passively tuned to be higher or lower than unity depending on the type of its nonlinearity (i.e., softening or hardening), with all the implications discussed above; specifically, the well-known capacity for higher harmonic generation and vibration energy scattering in the frequency domain directly affects the nonlinear bandwidth of nonlinear resonators, and, hence, their capacity to dissipate vibrations over different frequency bands. Hence, it is of great importance to study how the incorporation of an internal PID alters the bandwidth and the T-B product of an otherwise linear resonator (the PS). This has not been studied before and is the focus of this work.

Accordingly, the dissipative capacity of a PS incorporating an internal PID is studied herein employing the previous concepts of nonlinear bandwidth and T-B product. The manuscript is organized as follows: Section 2 introduces the integrated system of the PS incorporating a PID with different initial configurations of granules. Section 3 presents the computation of the nonlinear bandwidth of the integrated PS-PID system, as well as the corresponding T-B product. Section 4 focuses on various PID configurations possessing a fixed number of granules and discusses differences from the bandwidth perspective between these configurations. Section 5 estimates the time-bandwidth of the PS coupled to a PID of optimal configuration, whereas Section 6 presents some concluding remarks.

## 2. Discrete element modeling of the integrated PS-PID system

A single-degree-of-freedom (single-DOF) linear oscillator (namely, the primary structure – PS), incorporating an internal PID (essentially acting as the NES) with a varying number of granules (in different initial topologies) shown in their initial equilibrium configurations, are depicted in Fig. 1. The mass, stiffness, and damping of the PS are denoted by $m_{ps}$, $K$, and $C$, respectively, and it is subject to the shock force $F(t) = F_0 \sin\left[(\pi/t_0) \cdot \min(t, t_0)\right]$, where $F_0$ and $t_0$ are the shock amplitude and duration, respectively. A two-dimensional (2D) rectangular PID container with length $d_1$ and height $d_2$ is situated within the PS, cf. Fig. 1a, whose walls are assumed to be rigid; additional geometric parameters indicating the initial clearances between granules are shown in Fig. 1b-d; these are the same configurations introduced in [3]. The planar PIDs (the dynamics are assumed to be in a 2D plane) are constructed by placing various numbers of granules (particles) at different initial equilibrium positions or topologies inside the rectangular container. No pre-compression exists between the granules and the walls or between the granules themselves. Moreover, all granules are spherical and identical, are composed of the same linear viscoelastic material, e.g., steel, and are assumed to move in a 2D horizontal level plane; accordingly, the weights of the granules have no effect on the granular dynamics, so gravitational forces are neglected in this work.

Fig. 1b depicts three initial configurations for a PID with a single granule, namely, asymmetric middle (Config. 1), and two asymmetric (close to either the right or left wall, Configs. 2 and 3) positions. It is noted that only the initial position is changed in these three configurations to study the possible effects on the dynamics. The clearance between the granule and the wall in Config. 1 is denoted by $d_0$, whereas the clearance between the granule and the opposite wall in the other two configurations is equal to $2d_0$. Similarly, for a PID with two granules (cf. Fig. 1c), two symmetric (cf. Configs. 1 and 4) and four asymmetric initial configurations (cf. Configs. 2, 3, 5, and 6) are considered. Lastly, three initial configurations of PIDs with three, five and eight granules are presented in Fig. 1d; in this case the granules are arranged in topologies with clearance $d_0$ between neighboring granules and gap $d_v$ between the top granule and the ceiling of their container. In all cases the granules are identical and spherical, and composed of steel. After applying a shock excitation, the granules oscillate inside their 2D container, undergoing granule-wall and granule-granule collisions which are interrupted by periods of "free flight". Moreover, the displacement of the PS is denoted by the variable *z*.

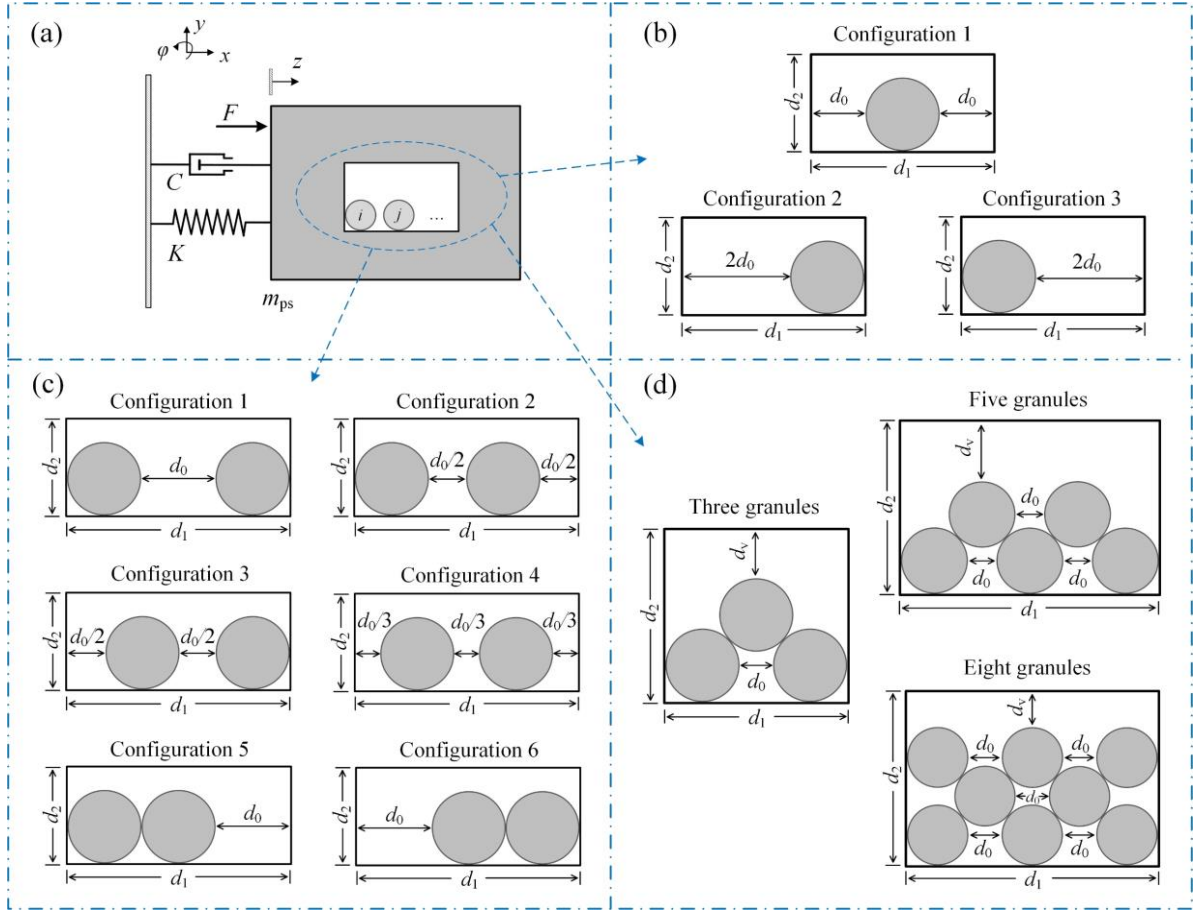

**Fig. 1.** Schematics of the considered nonlinear resonators: (a) Dissipative linear PS with a planar PID, subject to a shock force; (b) PID composed of a single granule in three distinct initial configurations; (c) PID composed of two granules in six initial configurations; and (d) PID composed of multiple granules in three initial configurations [3].

The integrated PS-PID system possesses highly nonlinear and discontinuous (non-smooth) dynamics caused by granule-wall and granule-granule collisions. To this end, the discrete element method (DEM) is a useful computational tool to accurately simulate the transient response following the applied shock; this is also a basic prerequisite for exploring the nonlinear bandwidth characteristics of this system. DEM is based on the Hertzian contact law [37] assuming small elastic deformations in granule-granule and granule-wall contact interactions within a small enough simulation time step. In each collision, there are strongly nonlinear forces developed, namely, a nonlinear (compressive Hertzian) contact force in the normal direction, and a frictional force in the tangential direction [3,38–40]. Therefore, due to these highly nonlinear and discontinuous contact forces, when the shock energy is applied to the PS, there can be nonlinear targeted energy transfer (TET) from the PS to the PID, whereby the input shock energy is irreversibly transferred and locally dissipated by both granular inelastic collisions and frictional effects. Accordingly, following the computational model developed in [3], the governing equations of motion of the integrated system are given by,

$$
\begin{aligned}
&m_{\text{PS}} \ddot{z} + C\dot{z} + Kz = F + F_{\text{d}} \\
&m_i \ddot{\mathbf{u}}_i = \sum_j \left( \mathbf{N}_{ij} + \mathbf{f}_{ij} \right) + \sum_k \left( \mathbf{N}_{ik} + \mathbf{f}_{ik} \right) \\
&I_i \ddot{\boldsymbol{\theta}}_i = R_i \sum_j \left( \mathbf{n}_{ij} \times \mathbf{f}_{ij} \right) + R_i \sum_k \left( \mathbf{n}_{ik} \times \mathbf{f}_{ik} \right)
\end{aligned} \quad (1)
$$

where overdot indicates derivation to time $t$; referring to Fig. 2, the subscripts $i$ and $j$ denote the $i$-th and $j$-th granule, respectively; the subscript $k$ represents the $k$-th point on the wall of the PID container; $\mathbf{u}_i$ and $\boldsymbol{\theta}_i$ characterize the $i$-th granule displacement vector and angular displacement pseudo-vector, respectively; $m_i$, $I_i$, and $R_i$ are the mass, moment of inertia and radius of granule $i$, respectively; $\mathbf{N}_{ij}$ and $\mathbf{N}_{ik}$ are the normal contact forces between granules $i$ and $j$, and particle $i$ and the point on the wall $k$, respectively; $\mathbf{f}_{ij}$ or $\mathbf{f}_{ik}$ represent the corresponding tangential forces; $\mathbf{n}_{ij}$ and $\mathbf{n}_{ik}$ are the unit vectors pointing from granule $i$ to granule $j$, and from granule $i$ to the point on the wall $k$, respectively; and lastly, $F_{\text{d}}$ is the magnitude of the resultant contact forces in $x$ direction acting on the PID container walls by all granules.

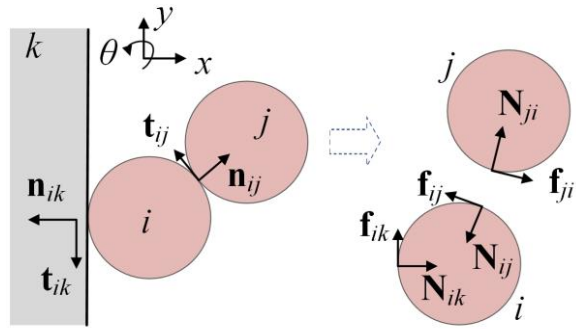

**Fig. 2.** Force analysis diagram for granule-granule and granule-wall contact interactions [3].

Adopting an inelastic Hertzian law for the radial contact forces $\mathbf{N}_{ij}$, and a continuous smooth Coulomb-tanh friction model for the tangential contact forces $\mathbf{f}_{ij}$ [3,14,39,41] to avert possible numerical instability, we set,

$$\mathbf{N}_{ij} = -\left(A_{ij}\delta_{n,ij}^{\frac{3}{2}} + \gamma_{ij}\dot{\delta}_{n,ij}\right)\mathbf{n}_{ij}$$
$$\mathbf{f}_{ij} = -\mu\left|\mathbf{N}_{ij}\right|\tanh\left(k_s\dot{\delta}_{t,ij}\right)\mathbf{t}_{ij} \tag{2}$$

where $A_{ij} = \frac{4}{3}E_{\text{eff}}\sqrt{R_{\text{eff}}}$ is the contact efficient; $\gamma_{ij} = \alpha_n\sqrt{m_{\text{eff}}A_{ij}}\delta_{n,ij}^{\frac{1}{4}}$ is the damping coefficient; $\delta_{n,ij} = \max\left(R_i + R_j - |\mathbf{u}_j - \mathbf{u}_i|, 0\right)$ is the penetration depth between granules $i$ and $j$; the unit vector $\mathbf{n}_{ij}$ is defined by $\mathbf{n}_{ij} = (\mathbf{u}_j - \mathbf{u}_i)/|\mathbf{u}_j - \mathbf{u}_i|$; $\mu$ is the slipping friction coefficient; $k_s$ controls the smoothness of the friction and also models the frictional forces close to zero velocity; $\mathbf{t}_{ij}$ is the tangential unit vector which is perpendicular to $\mathbf{n}_{ij}$; and $\dot{\delta}_{t,ij} = \left[\left(\dot{\mathbf{u}}_i + R_i\dot{\boldsymbol{\theta}}_i \times \mathbf{n}_{ij}\right) - \left(\dot{\mathbf{u}}_j + R_j\dot{\boldsymbol{\theta}}_j \times \mathbf{n}_{ji}\right)\right] \cdot \mathbf{t}_{ij}$ is the scalar relative velocity between granules $i$ and $j$. Moreover, we define the effective mass, radius, and Young's modulus, $m_{\text{eff}}$, $R_{\text{eff}}$, and $E_{\text{eff}}$, respectively, as $\frac{1}{E_{\text{eff}}} = \frac{1-v_i^2}{E_i} + \frac{1-v_j^2}{E_j}$, $\frac{1}{R_{\text{eff}}} = \frac{1}{R_i} + \frac{1}{R_j}$, and $\frac{1}{m_{\text{eff}}} = \frac{1}{m_i} + \frac{1}{m_j}$. Here, $E_{i(j)}$ and $v_{i(j)}$ represent the Young's modulus and Poisson's ratio of the granule $i$ ($j$), respectively. When considering a contact between granule $i$ and the point on the container wall $k$, this can be treated as contact between two identical granules. Therefore, the normal force $\mathbf{N}_{ik}$ and the tangential force $\mathbf{f}_{ik}$ can also be computed using Eq. (2) (see Refs. [3,39]).

For PIDs with a single or two granules (cf. Figs. 1b,c), the frictional force is neglected since the dynamics is one-dimensional (1D), whereas the slipping friction due to granular rotations is taken into account for the PIDs with more granules (cf. Fig. 1d). In the latter cases, the granules undergo translational motions in the $x$ and $y$ directions and rotations measured by the corresponding angles. The equations of motion (1) are numerically solved by a 4$^{\text{th}}$ order Runge-Kutta algorithm in MATLAB® based on the variable time step scheme developed in [3]. This is required to accurately capture the times of granular collisions, so a relatively large time step $\Delta t_1$ is adopted for "free flight," and a smaller time step $\Delta t_2$ is used for capturing collision interactions; in addition, a smooth time step switch is imposed between the large to small time steps. A criterion on energy convergence is followed to ensure an accurate DEM simulation process; that is, the summation of all energy components (i.e., instantaneous potential and kinetic and accumulated dissipative energies) should be conserved at each time instant. Lastly, to get reliable simulation results, the time step $\Delta t_2$ and the friction smoothing parameter $k_s$ – see (2), are selected to ensure convergent DEM simulation results.

The total mass of the integrated PS-PID system (including the cases with and without the PIDs) is kept fixed to $M$ to avoid possible mass-added effects. The system parameters of the PS without the PID are listed in Table 1, whereas the mass ratio of the granules relative to the mass of the PS is fixed to 6%. Some of the system parameters of the (identical) granules of the PID are listed in Table 2. The PS-PID system with zero initial conditions is excited by the applied shock, following which the nonlinear transient dynamics ensue. Unless otherwise noted, the shock magnitude $F_0$ and its duration $t_0$ are fixed to $F_0 = 5\times10^3$ N and $t_0 = 1\times10^{-3}$ s.

Table 1. The system parameters of the PS without the PID

| Item | Value |
| --- | --- |
| Mass, $M$ [kg] | 20 |
| Linear stiffness coefficient, $K$ [N·m$^{-1}$] | $8\times10^4$ |
| Linear damping coefficient, $C$ [N·s·m$^{-1}$] | 25.30 |

Table 2. The system parameters of the granules of the PID

| Item | Value | Parameter | Value |
| --- | --- | --- | --- |
| Total mass ratio, $\varepsilon$ | 6% | Poisson's ratio, $v$ | 0.3 |
| Young's modulus, $E$ [Pa] | $200\times10^9$ | Coefficient of restitution, $\alpha_n$ | $6.313\times10^{-3}$ |
| Granule's density, $\rho$ [kg·m$^{-3}$] | 7,850 | Frictional coefficient, $\mu$ | 0.099 |

## 3. Computation of the nonlinear bandwidth

As mentioned previously, the bandwidth of the decaying response of the integrated PS-PID system is directly related to the rate of energy dissipation of this system, which is expected to be energy-dependent due to the strong nonlinearity of the transient dynamics. According to [35,36], the nonlinear bandwidth is computed based on the instantaneous mechanical energy of the system; moreover, as discussed previously, the various energy components are also needed to monitor energy conservation (and thus accurate convergence) of the numerical algorithm.

The mechanical energy of the PID is denoted by $E_{\mathrm{PID}}$, and is computed by adding the instantaneous energies of all granules, $E_{\mathrm{p},i}$. In addition, the total instantaneous energy of the PS is denoted by $E_{\mathrm{PS}}$. The energy expressions are given by

$$E_{\mathrm{PS}}(t) = \frac{1}{2} m_{\mathrm{PS}} \left(\dot{z}(t)\right)^2 + \frac{1}{2} K \left(z(t)\right)^2,$$
$$E_{\mathrm{PID}}(t) = \sum_i E_{\mathrm{p},i}(t), \qquad (3)$$
$$E_{\mathrm{p},i}(t) = \frac{1}{2} m_i |\dot{\mathbf{u}}_i(t)|^2 + \frac{1}{2} I_i \left(\dot{\theta}_i(t)\right)^2 + \sum_j \frac{1}{5} A_{ij} \left(\delta_{\mathrm{n},ij}(t)\right)^{\frac{5}{2}} + \sum_k \frac{2}{5} A_{ik} \left(\delta_{\mathrm{n},ij}(t)\right)^{\frac{5}{2}}$$

These instantaneous energy measures can be normalized with respect to the total input shock energy, $E_{\mathrm{in}} = \int_0^{t_0} F(\tau)\dot{z}(\tau)\mathrm{d}\tau$. The corresponding (percentages of) normalized instantaneous energies of the PID, PS, and integrated PS-PID system are denoted by $\eta_{\mathrm{PID}}$, $\eta_{\mathrm{PS}}$, and $\eta_{\mathrm{sys}}$, respectively, as

$$\eta_{\mathrm{PID}}(t) = \frac{E_{\mathrm{PID}}(t)}{E_{\mathrm{in}}} \times 100\%, \quad \eta_{\mathrm{PS}}(t) = \frac{E_{\mathrm{PS}}(t)}{E_{\mathrm{in}}} \times 100\%, \quad \eta_{\mathrm{sys}}(t) = \frac{E_{\mathrm{PID}}(t) + E_{\mathrm{PS}}(t)}{E_{\mathrm{in}}} \times 100\% \qquad (4)$$

Regarding energy dissipation, the shock energy is dissipated partially by the damping of the PS, and by viscoelastic and frictional forces during granule-granule or granule-wall collisions. Therefore, shock energy transferred from the PS to the PID is locally dissipated by the PID owing to inelastic granular collisions and slipping friction due to relative granular rotations; these components of dissipative energy up to time instant $t$ are represented by $W_{\mathrm{vis}}(t)$ and $W_{\mathrm{f}}(t)$, respectively. Lastly, the cumulative dissipative energies up to time instant $t$ by the PID and PS are denoted by $W_{\mathrm{PID}}(t)$ and $W_{\mathrm{PS}}(t)$, respectively, as

$$W_{\text{PS}}(t) = \int_0^t C\left(\dot{z}(\tau)\right)^2 \mathrm{d}\tau,$$

$$W_{\text{PID}}(t) \equiv W_{\text{vis}}(t) + W_{\text{f}}(t),$$

$$W_{\text{vis}}(t) = -\int_0^t \left[\sum_i \left(\sum_j \gamma_{ij} \dot{\delta}_{\text{n},ij}(\tau)\mathbf{n}_{ij} \cdot \dot{\mathbf{u}}_i(\tau) + \sum_k \gamma_{ik} \dot{\delta}_{\text{n},ik}(\tau)\mathbf{n}_{ik} \cdot \dot{\mathbf{u}}_i(\tau)\right)\right]\mathrm{d}\tau, \quad (5)$$

$$W_{\text{f}}(t) = -\int_0^t \left[\sum_i \left(\sum_j \mathbf{f}_{ij}(\tau) \cdot \left(\dot{\mathbf{u}}_i(\tau) + R_i\dot{\theta}_i(\tau)\mathbf{t}_{ij}\right) + \sum_k \mathbf{f}_{ik}(\tau) \cdot \left(\dot{\mathbf{u}}_i(\tau) + R_i\dot{\theta}_i(\tau)\mathbf{t}_{ik}\right)\right)\right]\mathrm{d}\tau$$

Hence, the normalized cumulative energy dissipation measures for the PID and PS are

$$r_{\text{PID}}(t) = \frac{W_{\text{PID}}(t)}{E_{\text{in}}} \times 100\%, \quad r_{\text{PS}}(t) = \frac{W_{\text{PS}}(t)}{E_{\text{in}}} \times 100\% \quad (6)$$

When the time window [0, *t*] is large enough, the measures above are referred to as *eventual normalized dissipation measures* since they characterize the dissipation over the entire duration of the transient dynamics.

The nonlinear bandwidth measure is now considered to characterize the dissipative capability of the integrated PS-PID system. The nonlinear bandwidth or RMS bandwidth denoted by $\Delta\omega^{**}$ is based on the variance of the mechanical energy in the frequency domain and is expressed as [35,36],

$$\Delta\omega^{**} = 2\sqrt{\frac{\int_{-\infty}^{\infty} \omega^2 E^2(\omega)\mathrm{d}\omega}{\int_{-\infty}^{\infty} E^2(\omega)\mathrm{d}\omega}} \quad (7)$$

where $E(\omega)$ is the energy spectrum of the system at frequency $\omega$. Specifically, considering the velocity time series of the PS, namely $\dot{z}$, the energy spectrum is defined as $E(\omega) = \left|\mathbb{F}\left[\langle \dot{z} \rangle\right]\right|^2$, where $\langle \cdot \rangle$ denotes the envelope operator, and $\mathbb{F}$ the Fourier transform operator. Note that the nonlinear bandwidth computation may be based on the decaying energy of the PS alone, or, alternatively, on the decaying energy envelope of the overall integrated PS-PID system. Typically, the (decaying) envelope of the velocity of the PS is employed because its envelope coincides with the corresponding PS energy decay [35,42]; however, in the following analysis the nonlinear bandwidth of the integrated system will be considered as well for comparison. In addition, when computing the nonlinear bandwidth, the prerequisite condition is that the energy remaining in the PS should approach zero asymptotically, i.e., that the PS undergoes a *decaying transient oscillation*.

Following the bandwidth computation, the time-bandwidth (T-B) product for the PS can be estimated. As mentioned previously, the T-B product for a classical linear, damped resonator is always equal to unity, i.e., $\Delta\tau\,\Delta\omega = 1$, where $\Delta\tau$ is the decay-time constant (storage time), and in this case $\Delta\omega$ is the classical half-power bandwidth. The T-B product of the strongly nonlinear dynamics of the PS incorporating a PID is also computed as $\Delta\tau^{**}\Delta\omega^{**}$ where (7) is considered, and $\Delta\tau^{**}$ is the time required for the energy to decay to $(1/e)$ of its initial value. As mentioned in the Introduction, for passive nonlinear time-invariant oscillators the previous classical limit of unity can be broken either above or below, with important implications on the capacity of the nonlinear oscillator to store or dissipate mechanical energy over certain frequency bands. In the following exposition a detailed study of the nonlinear bandwidth and the T-B product of the PS incorporating a PID with varying number of granules and different initial topologies will be studied. To the authors' knowledge this is the first such study appearing in the literature in that context, and the outlined methodology and obtained results can be employed to study the capacity for energy dissipation by the PIDs considered.

## 4. Nonlinear bandwidth and T-B product for PIDs with one or two granules

Considering initially the PIDs consisting of either a single granule or two granules, we study the initial configurations depicted in Figs. 1b,c, respectively. The dynamics in these cases are simpler as the granules move only in the *x* direction on a horizontal plane, so granular rotations and frictional forces are omitted. The corresponding system parameters are listed in Tables 1-3. In [3] the convergence of the computational algorithm and energy conservation were studied, yielding the large time step $\Delta t_1 = 2\times 10^{-5}$ s and the small time step $\Delta t_2 = 3\times 10^{-8}$ s for the DEM simulations.

In Fig. 3 the nonlinear bandwidth, decay-time constant (i.e., storage time) and the resulting T-B product for the PS incorporating a PID with single granule, are depicted as functions of the input shock energy (corresponding to shock magnitudes $F_0 = 10$ N, 50 N, 100 N, 500 N, $1\times 10^3$ N, $2.5\times 10^3$ N, $5\times 10^3$ N, $7.5\times 10^3$ N and $1\times 10^4$ N). These measures were computed by Eq. (7) based on the DEM simulation results of the decaying response of the PS, and their strong dependence (tunability) on the input energy is evident. The performance of the three single-granule PID initial configurations (cf. Fig. 1b) are compared to the linear case of the PS system without PID (but with the same total mass). Considering the linear case as reference baseline, the bandwidth, decay-time constant and T-B product are equal to 1.26 rad/s, 0.79 s and unity, respectively.

Table 3. System parameters of a PID with a single or two granules.

| Parameter | PID with a single granule | PID with two granules |
|---|---|---|
| Granule's mass, $m$ [kg] | 1.2 | 0.6 |
| Granule's radius, $R$ [mm] | 33.2 | 26.3 |
| Clearance, $d_0$ [mm] | 0.4 | 6.0 |
| Length, $d_1$ [mm] | 67.2 | 111.2 |

For small enough input energy ($F_0 < 500$ N) the time-bandwidth measures for initial Config. 1 are identical to the linear case due to the absence of any granule-wall collisions at this low energy level. As the input energy is increased, however, the time-bandwidth measures show strong tunability with energy, but also (perhaps surprisingly) high sensitivity to the initial granule configurations. Indeed, at relatively small shock excitations, the bandwidth and decay-time constant are different for the three initial configurations, but as the shock magnitudes are increased, these measures start converging.

For example, at the low shock magnitude $F_0 = 100$ N (highlighted in Fig. 3), the nonlinear bandwidth for Config. 3 is much smaller compared to the other two configurations, which might indicate that the dissipation capacity of the PID with initial configuration 3 is diminished compared to the other two initial configurations, even though the initial energy of the PS drops to (1/e) of its value rather quickly – judging by Fig. 3b. It is interesting to note, however, that for relatively strong shocks the results converge for the three configurations, so the initial granule configuration plays a near negligible role at these high input energies, as expected by physical intuition. Lastly, at low energies we note that the T-B product of the PS for Configs. 2 and 3 decreases significantly below the classical linear limit of unity (as the decay-time constant also decreases for these cases); this shows diminished overall dissipative capacity of the PID at low shocks. In other words, even though the initial energy of the PS drops significantly, as evidenced by the smaller values of $\Delta\tau^{**}$ compared to the linear PID-less case), the oscillations of the PS persist for longer time (compared to the linear case) even though they exhibit small amplitudes.

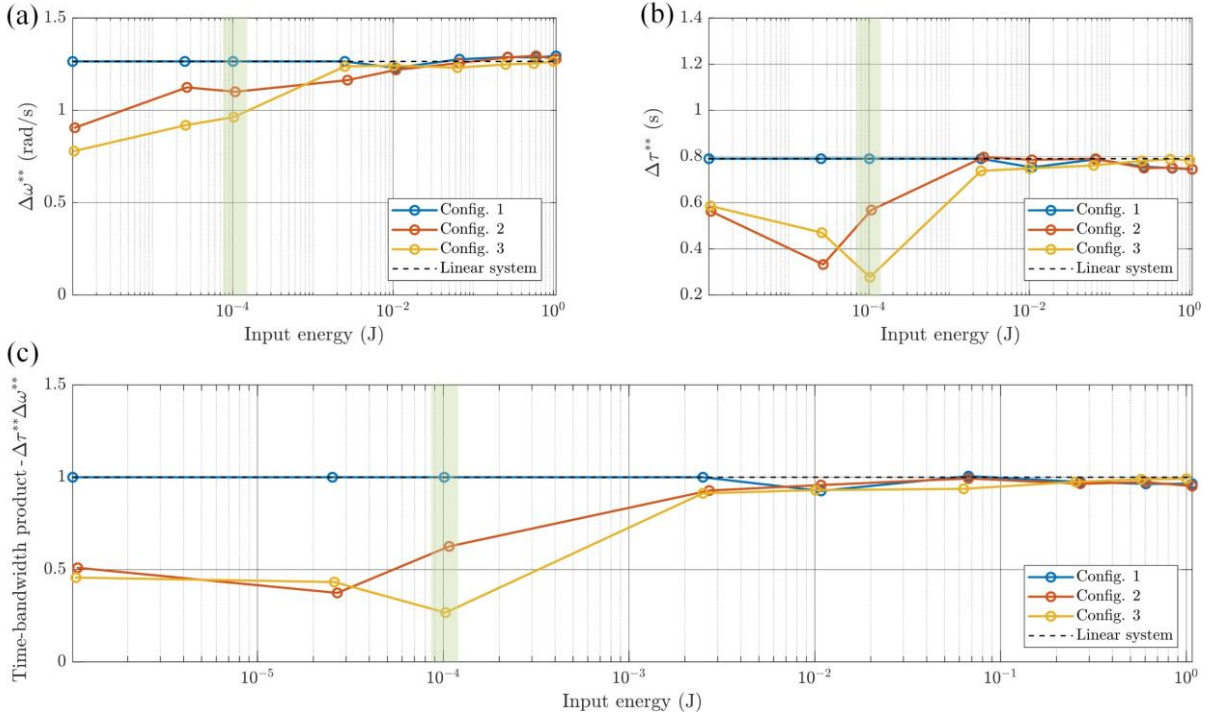

**Fig. 3.** PS with a PID possessing a single granule at different initial configurations, cf. Fig.1b: (a) Nonlinear bandwidth, (b) decay-time constant, (c) T-B product for varying input shock energy; these computations are based on the decaying response of the PS (the shock force magnitude corresponding to $F_0 = 100$ N is marked by the highlighted area, cf. Fig.4).

The results of Fig. 3 show that, compared to the linear case (e.g., with no PID or low-energy Config. 1), the PS incorporating a PID with initial Config. 2 stores vibrations for longer times (i.e., it possesses weaker overall dissipative capacity) and with smaller bandwidths. In fact, as discussed below (also cf. Fig. 4), at early times the motion of the PS decays faster compared to the linear case (PS with no PID), but eventually the overall decay rate of the integrated PS-PID system is slower compared to the linear decay. This is because in the early stage of the transient dynamics a part of the input energy is stored (and preserved for short time periods) as kinetic energy in the granule during its phase of "free flight." This energy is eventually transferred back to the PS and dissipated at later times, but this process delays the overall decay rate of the PS oscillation and yields a weaker overall dissipative capacity for the system. This effect is encountered in PIDs with one or two granules where the possibility of "free flight" exists, but to a lesser extent in PIDs with a higher number of granules where such granular motions are more restricted.

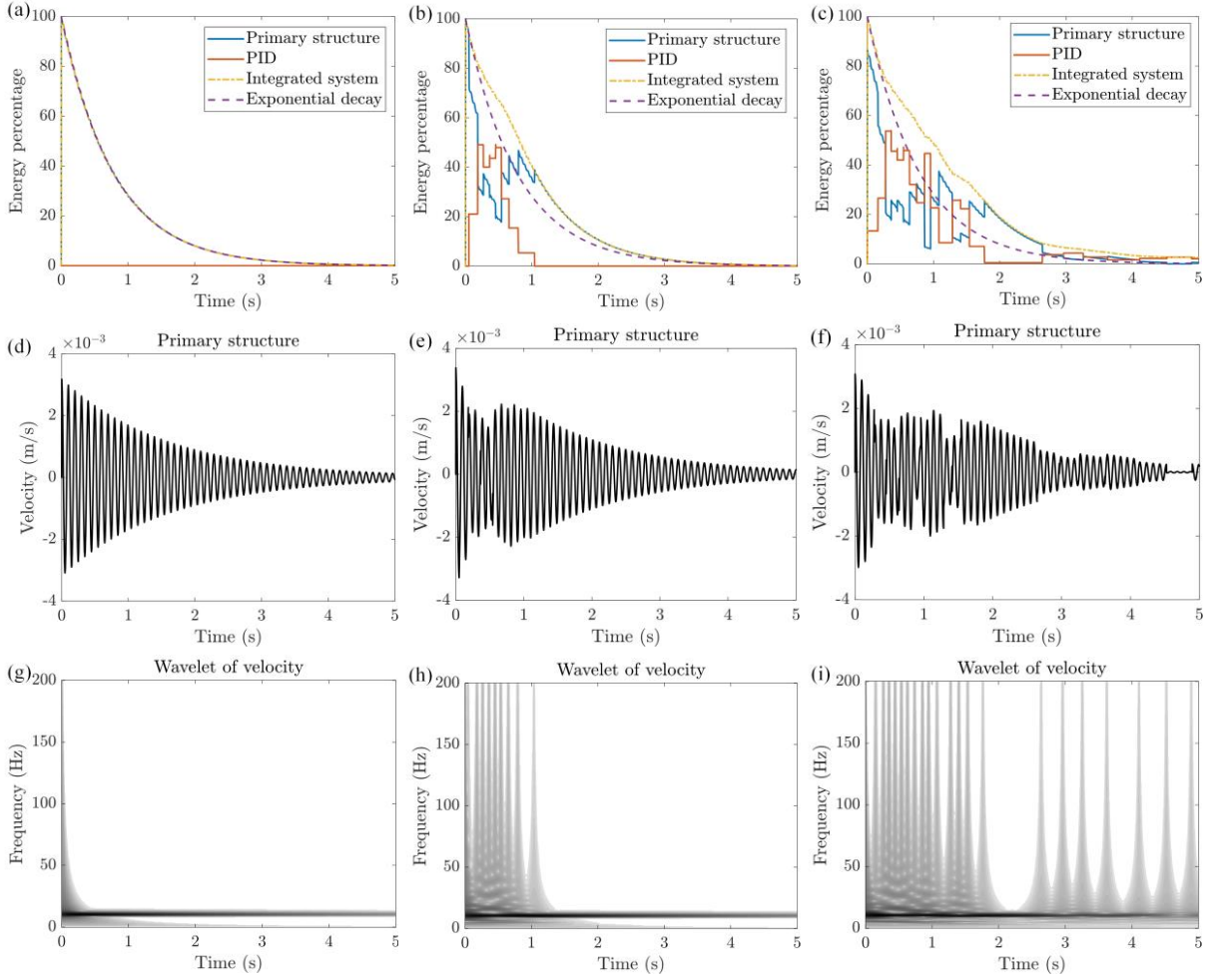

**Fig. 4.** PS incorporating a PID with a single granule for $F_0 = 100$ N: (a-c) Instantaneous normalized energy percentage, (d-f) velocity of the PS, and (g-i) its normalized wavelet transform spectrum; (a), (d) and (g) Config. 1; (b), (e) and (h) Config. 2; and (c), (f) and (i) Config. 3.

Focusing on the low-intensity shock amplitude at $F_0 = 100$ N, the corresponding normalized energy percentages for the three initial configurations are depicted in Figs. 4a,b, and c, respectively. For comparison the exponential decay of energy (denoted by the dashed curve) corresponding to the PS without a PID is also shown. Fig. 4a shows that no PS-PID energy exchanges occur for Config. 1 due to the absence of granule-wall collisions at this low shock level. The results of Figs. 4b,c reveal the much enhanced dissipation rate of the vibrations of the PS at the early stage of the motion compared to the linear case, which is due to the intense inelastic granule-wall collisions. At the same time, the decay rate of the integrated PS-PID system is smaller compared to linear exponential decay, since between granular-wall collisions there are periods of "free flight" of the granule where a part of the input energy is stored and preserved as kinetic energy in the granule instead of being dissipated. Eventually, this kinetic energy component of the granule motion is transferred back to the PS and is dissipated, but this delays the overall decay rate of the PS oscillation. The velocity $\dot{z}(t)$ of the PS and its corresponding normalized wavelet transform spectrum are depicted in Figs. 4d–i for all three initial configurations. Heavy (lighter) shades of the wavelet spectra indicate higher (smaller) energies, and the sudden high-frequency energy bursts indicate the intense input energy scattering from low-to-high frequencies that occur during the inelastic granule-wall collisions. As discussed in [35,36] such nonlinear energy transfers in the frequency domain directly affect the dissipative capacity of the system, and hence its nonlinear bandwidth.

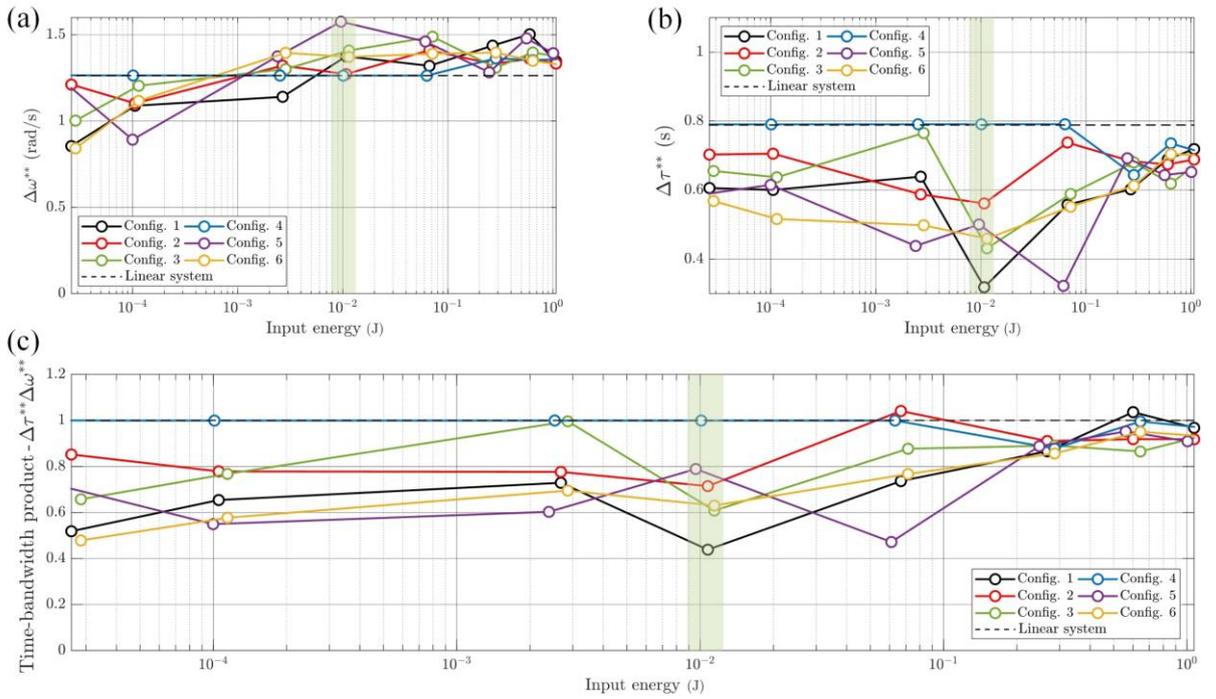

**Fig. 5.** PS with a PID possessing two granules at different initial configurations, cf. Fig.1c: (a) Nonlinear bandwidth, (b) decay-time constant, (c) T-B product for varying input shock energy; these computations are based on the decaying response of the PS (the shock force magnitude corresponding to $F_0 = 1 \times 10^3$ N is marked by the highlighted area, cf. Fig.6).

Next, the case of the PID with two granules (cf. Fig. 1c) is studied with system parameters listed in Table 3. There are six initial configurations shown in Fig. 1c, with granules at symmetric (Configs. 1 and 4) or asymmetric positions (Configs. 2, 3, 5, and 6). The corresponding time-bandwidth results for the PS as functions of input energy are depicted in Fig. 5. In terms of the nonlinear bandwidth, for relatively small shock amplitudes, e.g., at $F_0 = 1\times10^3$ N (marked by the highlighted areas in Fig. 5), substantial differences between the six initial configurations are observed. In particular, the bandwidth of the PS for Config. 2 is rather small among these cases and is close to the classical linear half-power bandwidth of the PS without a PID. Some differences are also found at relatively low shock amplitudes for the decay-time constant and the time-bandwidth product among the six configurations. The reason is that at relatively low shock levels the relatively *fast* granular motions dominate the system dynamics, compared to the relatively *slow* dynamics of the PS. Typically the T-B product becomes less than unity irrespective of the level of input energy, except for Configs. 1 and 2 that can become greater than unity at certain input energy levels. Also, the T-B product for Config. 4 is equal to unity at low input energies due to the absence of granule-wall interactions up to shock magnitudes greater than $2.5\times10^3$ N.

In Fig. 6 the normalized energy measures, and the velocity of the PS together with its normalized wavelet transform spectrum are presented for Configs. 1, 2, and 5 at the small shock magnitude $F_0 = 1\times10^3$ N. According to the depicted normalized energy measures of Figs. 6a,b,c, the energy of the PS decays faster compared to the linear case, yielding a stronger dissipative capacity of the PID in these cases. In Figs. 6g,h,i, intense input energy scattering from low to high frequencies is realized, as evidenced by the "energy bursts" at each granular collision, especially in the early time, highly energetic regime of the transient dynamics. This high frequency energy scattering is much more intense compared to the PID with a single granule (cf. Fig. 4) and explains the stronger dissipative capacity of the two-granule PID (since by transferring a significant portion of the input energy at higher frequencies, one achieves lower vibration amplitude and much more enhanced energy dissipation).

Based on these results it is reasonable that the nonlinear bandwidth of the PS for each of these configurations (cf. Fig. 5a) is larger compared to the linear PID-less case. Furthermore, for Configs. 2 and 4 the nonlinear bandwidth of the PS is nearly identical to the linearized bandwidth and much smaller compared to the other configurations; this can be explained by the decay of the normalized energy of the PS in Fig. 6b which (after an early time regime where, as for the single granule – PID, kinetic energy in the granules is stored and preserved, with few collision interactions) is close to the linear exponential decay rate and the overall decay rate of the integrated PS-PID system. Therefore, the nonlinear bandwidths of the PS for Configs. 2 and 4 at this input energy level are relatively low, and the PID is not highly efficient in terms of dissipative capacity.

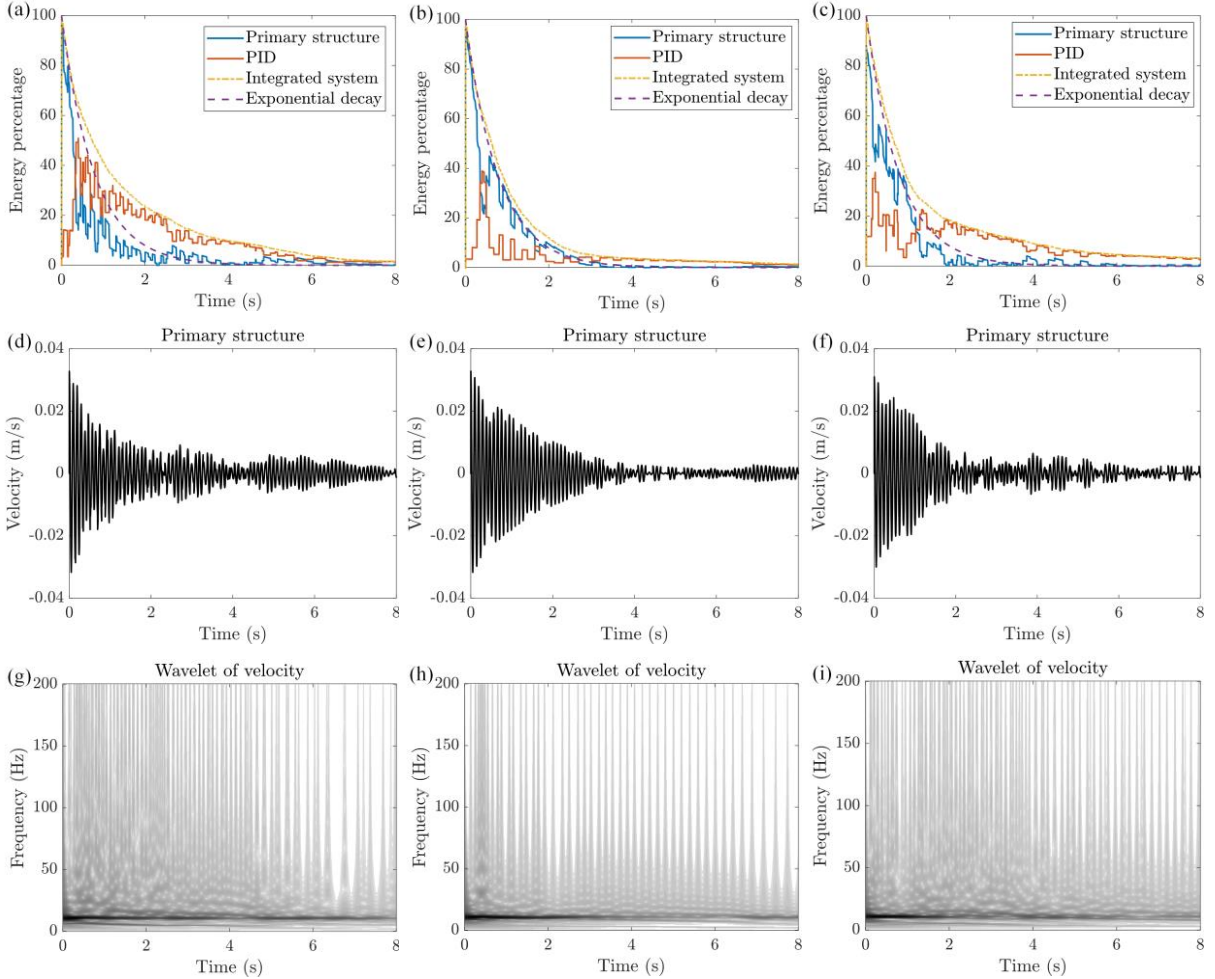

**Fig. 6.** PS incorporating a PID with two granules for $F_0 = 1\times10^3$ N: (a-c) Instantaneous normalized energy percentage, (d-f) velocity of the PS, and (g-i) its normalized wavelet transform spectrum; (a), (d) and (g) Config. 1; (b), (e) and (h) Config. 2; and (c), (f) and (i) Config. 5.

**5. Nonlinear bandwidth and T-B product for PIDs with more than two granules**

Next, we study cases of PIDs with more than two granules. The initial configurations of these PIDs, shown in Fig. 1d, yield dynamics that are much more complex [3], since they are 2D and involve additional frictional forces due to relative granular rotations (which were neglected in the 1D dynamics of the PIDs studied in Section 4). In the considered topologies the clearance between each granule and its neighbors is denoted by $d_0$, whereas the smallest gap between the granules and the ceiling of the container is given by $d_v$. The length and height of the PID container, namely $d_1$ and $d_2$, may vary according to the changes in $d_0$ and $d_v$, whereas keeping the initial topology unchanged. The granules move in both $x$ and $y$ directions and may undergo rotations yielding frictional forces during granule-granule or granule-wall interactions. Again, following [3] where detailed convergence studies of the numerical simulations were reported, the smoothing parameter in the tanh-friction law is taken as $k_s = 250$ s/m, the small time step as $\Delta t_2 = 3 \times 10^{-8}$ s, and the large time step as $\Delta t_1 = 2 \times 10^{-5}$ s.

Clearly, the size of the PID container as defined by the clearances $d_0$ and $d_v$ have major effect on the granular dynamics and, in turn, on the dissipative performance of the considered PIDs. With the aim of achieving maximum percentage of cumulative dissipated energy by the PID as defined by Eq. (6), a parametric study was performed in [3] to optimize the geometric parameters of the PID container for different numbers of granules. The main aim herein is to study time-bandwidth features corresponding to the optimal cases of the PIDs with multiple granules. This exercise was performed in [3] for the fixed shock magnitude $F_0 = 5 \times 10^3$ N, and the corresponding optimal cases are listed in Table 4; the optimal percentage of cumulative dissipated energy is 56.81% for three granules, 68.47% for five granules, and 70.85% for eight granules. Time-bandwidth computations of the optimal cases listed in Table 4 will be performed for different shock magnitudes, namely $F_0 = 100$ N, 500 N, $1 \times 10^3$ N, $2.5 \times 10^3$ N, $5 \times 10^3$ N, $7.5 \times 10^3$ N and $1 \times 10^4$ N.

Table 4. Optimal sizes of PIDs with multiple granules at $F_0 = 5 \times 10^3$ N [3]

| Parameters | Three granules | Five granules | Eight granules |
|---|---|---|---|
| Granule's mass, $m$ [kg] | 0.4 | 0.24 | 0.15 |
| Granule's radius, $R$ [mm] | 23.0 | 19.4 | 16.6 |
| Clearance, $d_0$ [mm] | 7.0 | 7.0 | 5.0 |
| Gap, $d_v$ [mm] | 11.4 | 29.9 | 22.6 |
| Length, $d_1$ [mm] | 99.0 | 130.4 | 109.5 |
| Height, $d_2$ [mm] | 95.0 | 100 | 110.0 |

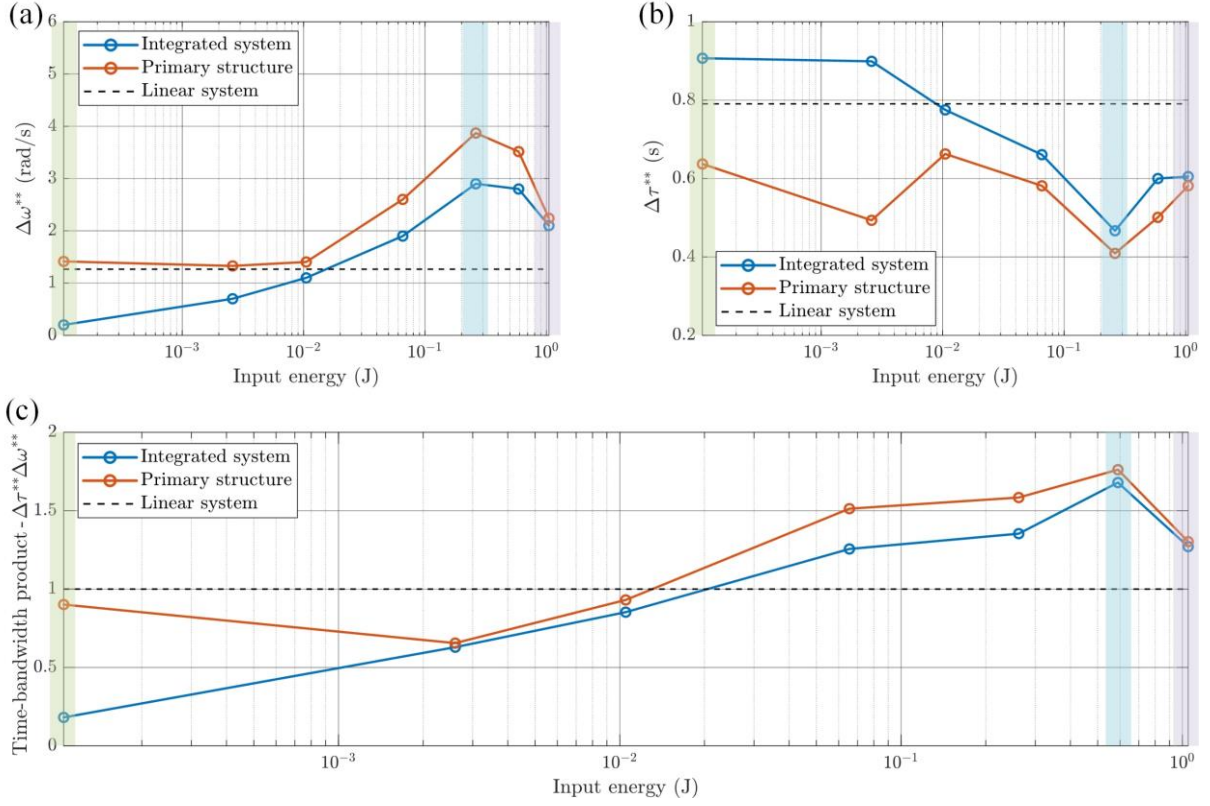

**Fig. 7.** PS with a PID possessing three granules for different shock excitations: (a) Nonlinear bandwidth, (b) decay-time constant, and (c) T-B product; the cases for shock magnitudes $F_0 = 100$ N, $F_0 = 5\times10^3$ N, and $F_0 = 1\times10^4$ N are marked by the green, blue, and purple shaded areas, respectively.

For the case of the optimal PID with three granules the nonlinear bandwidth, decay-time constant and T-B product are depicted in Fig. 7 as functions of the input shock energy (for force magnitudes ranging from 100 N to $1\times10^4$ N). These plots are based on the time series of the PS and the integrated PS-PID system, marked by orange and blue lines, respectively. At low shock excitations, the nonlinear bandwidth of the integrated system is smaller than that of the linear case, whereas the bandwidth of the PS is slightly larger than that of the linear case, cf. Fig. 7a. These results indicate that for weak applied shocks the PS dissipates shock energy faster than the linear case, but the converse holds for the integrated system. Furthermore, the corresponding T-B product for either the integrated system or the PS is below unity (e.g., at $F_0 = 100$ N marked by the green shaded area), the reason being that the decay-time constant of the PS is small, but the bandwidth of the integrated system is even smaller, eventually yielding T-B products less than the classical limit of unity. These results indicate that at low shocks, both the PS and integrated systems have overall capacities for storing vibrations for shorter times (i.e., they possess stronger overall dissipative capacity) *and* smaller bandwidths compared to the linear, PID-less case.

However, a different qualitative picture is obtained at more intense applied shocks, revealing the strong dependence (tunability) of the system on energy. Indeed, as the shock magnitude is increased, the T-B products for both the PS and the integrated system become greater than unity (cf. Fig. 7c). From the plots of Fig. 7b we note that the decay-time constant corresponding to the integrated system decreases with shock magnitude, which might mean stronger dissipative capacity; however, it turns out that this decrease is not in proportion with the corresponding increase of the nonlinear bandwidth (cf. Fig. 7a). The result is that, as the shock intensity is increased, the T-B product of the integrated system becomes greater than unity (principally due to the increase in the bandwidth), cf. Fig.7c, and the same holds for the PS response. This indicates that for strong shocks, both the PS and the integrated system attain overall capacities for storing vibrations for longer times (i.e., they possess weaker overall dissipative capacity) *and* larger bandwidths compared to the linear, PID-less case with the same bandwidth as the PS or the integrated system.

Hence, our quantitative time-bandwidth analysis emphasizes the drastic change in the dissipative performance of the three-granule PID with increasing input shock energy. It is interesting to note that the largest bandwidths, lowest decay-time constants and largest T-B products (greater than unity) for both the integrated system and PS are attained at $F_0 = 5\times10^3$ N (marked by the blue shaded area in Fig. 7); this coincides with the optimal energy dissipation performance for this system as reported in [3].

Next, the normalized energy measures and the velocity and its normalized wavelet transform spectrum are presented in Fig. 8 for the three-granule PID subject to force magnitudes of $F_0 = 100$ N (cf. Figs. 8a,d,g), $F_0 = 5\times10^3$ N (cf. Figs. 8b,e,h), and $F_0 = 1\times10^4$ N (cf. Figs. 8c,f,i). Note that the time-bandwidth results at these force magnitudes are highlighted by the shaded areas in Fig. 7. At the low shock $F_0 = 100$ N, a portion of the input energy is transferred into the granules and stored in the form of kinetic energy, cf. Fig. 8a, which is eventually dissipated by inelastic granular collisions and slipping friction. As a result, the PS dissipates shock energy rapidly, while the integrated system exhibits a slower dissipation rate. Additionally, it is observed from Fig. 8g that less input energy is being scattered in the high-frequency domain in this case. Therefore, the nonlinear time-bandwidth results for this low shock are verified, namely that the dissipative capacity of the integrated system is inferior compared to the linear case, whereas that of the PS is slightly enhanced, cf. Fig. 7a. As the shock amplitude or the input energy is increased, e.g., for $F_0 = 5\times10^3$ N and $1\times10^4$ N, intense energy exchanges between the PS and the PID occur (see Figs. 8b,c) and significant low-to-high energy scattering is realized (see Figs. 8h,i). Therefore, it is reasonable that at these higher shocks the corresponding nonlinear bandwidth is much larger than the linear case, and the storage time is smaller. Furthermore, there is intense input energy scattering from low-to-high frequencies (cf. Fig. 8h), further explaining the increased nonlinear bandwidth in these cases.

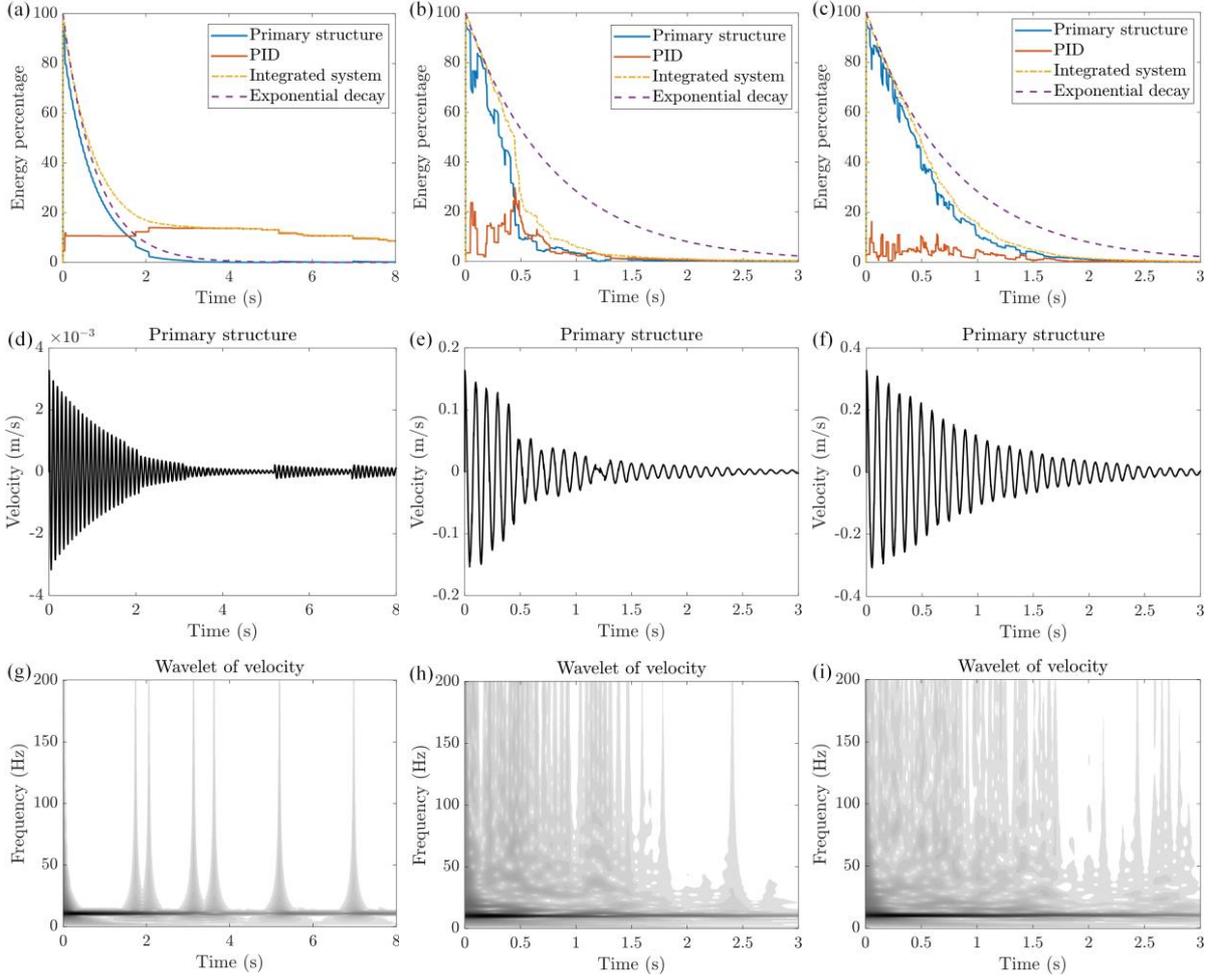

**Fig. 8.** PS incorporating a PID with three granules: (a-c) Instantaneous normalized energy percentage, (d-f) velocity of the PS, and (g-i) its normalized wavelet transform spectrum; (a), (d) and (g) $F_0 = 100$ N; (b), (e) and (h) $F_0 = 5\times10^3$ N; and (c), (f) and (i) $F_0 = 1\times10^4$ N.

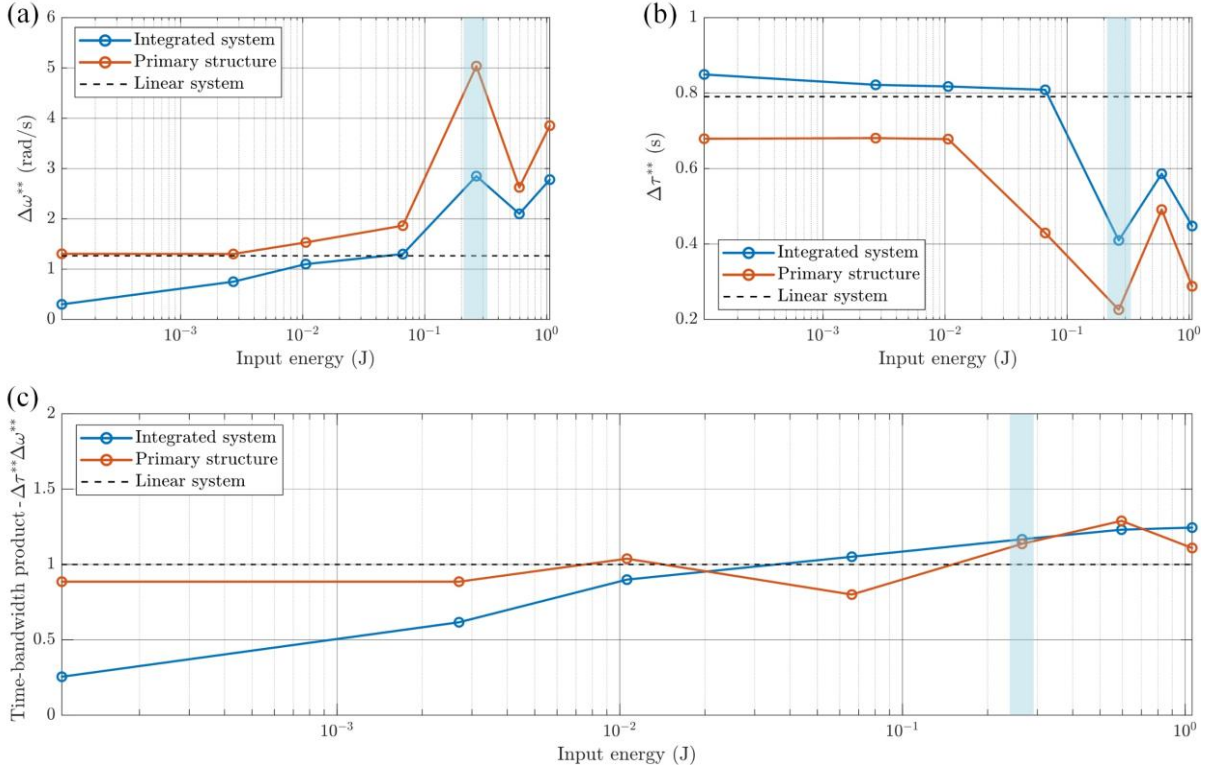

**Fig. 9.** Time-bandwidth for the case of five granules: (a) Nonlinear bandwidth, (b) decay-time constant, (c) T-B product; the shock force amplitude $F_0 = 5 \times 10^3$ N is marked by the blue shaded area.

Finally, similar trends are observed for the time-bandwidth results for the optimal cases of PIDs with five and eight granules, cf. Figs. 9 and 10, respectively. For example, at low shock amplitudes, the nonlinear bandwidth of the PS is slightly larger than in the linear system, but that of the integrated system is smaller compared to the linear system. As the intensity of the shock force is increased and the nonlinear effects of the granular dynamics become more pronounced, bandwidths of both the PS and the integrated system increase, becoming larger than the bandwidth predicted for the linear case. Therefore, the optimal PID configurations can attenuate the shock-induced vibration of the PS over a broad energy range, but the integrated PS-PID system dissipates shock energy less effectively at low shocks. In other words, for low shocks the integrated system can store the shock energy for longer time compared to the linear system. As the shock intensity increases, the decrease in the decay-time constant is not proportional to the increase of the nonlinear bandwidth, cf. Figs. 9 and 10, so the corresponding T-B product becomes greater than unity; this diminishes the overall dissipative capacity of the system compared to the linear system with no PID with the same bandwidth as the system with PID. Again, these results highlight the strong tunability of the dissipative capacity of the system with energy.

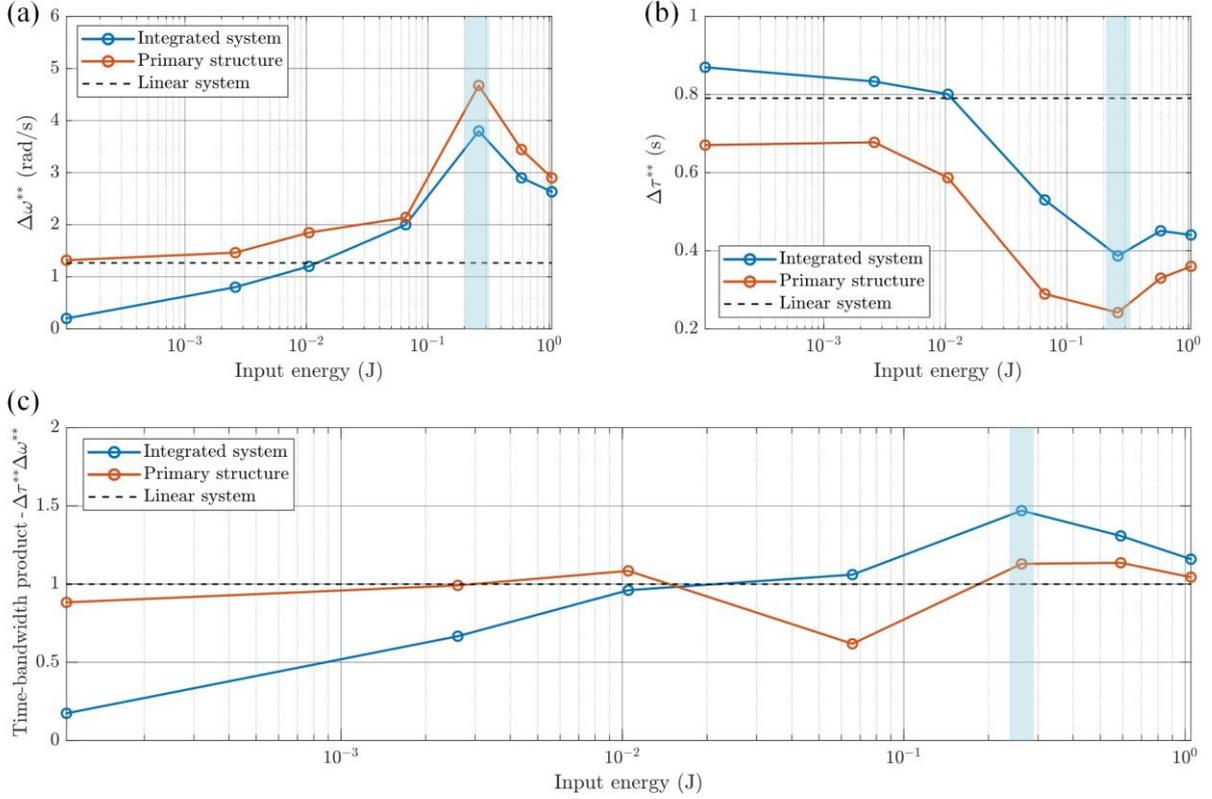

**Fig. 10.** Time-bandwidth for the case of eight granules: (a) Nonlinear bandwidth, (b) decay-time constant, (c) T-B product; the shock force amplitude $F_0 = 5\times10^3$ N is marked by the blue shaded area.

In Figs. 9 and 10, the bandwidths for the optimal cases of five or eight granules at $F_0 = 5\times10^3$ N (marked by the blue shaded areas) are the largest among these shock amplitudes. We note that the nonlinear bandwidths of the PS for the optimal cases of three, five and eight granules at $F_0 = 5\times10^3$ N are 3.87, 5.04 and 4.67 rad/s, respectively. Therefore, the nonlinear bandwidth of the optimal case of five granules at $F_0 = 5\times10^3$ N is the largest compared to the other optimal cases, which can be explained from the energy plots and wavelet transform spectra of Fig. 11. Indeed, there is irreversible energy transfer of input shock energy from the PS to the PID (cf. Figs. 11a,d) followed by local dissipation through inelastic granular collisions and friction effects. Furthermore, intense shock energy scattering from low-to-high frequencies is observed in the optimal cases of five or eight granules (cf. Figs. 11c,f) compared to the optimal case of three granules at $F_0 = 5\times10^3$ N (see Fig. 8h). Therefore, among the optimal cases listed in Table 4, the optimal PIDs with five and eight granules have more effective shock mitigation capacity compared to the optimal PID with three granules (cf. [3]).

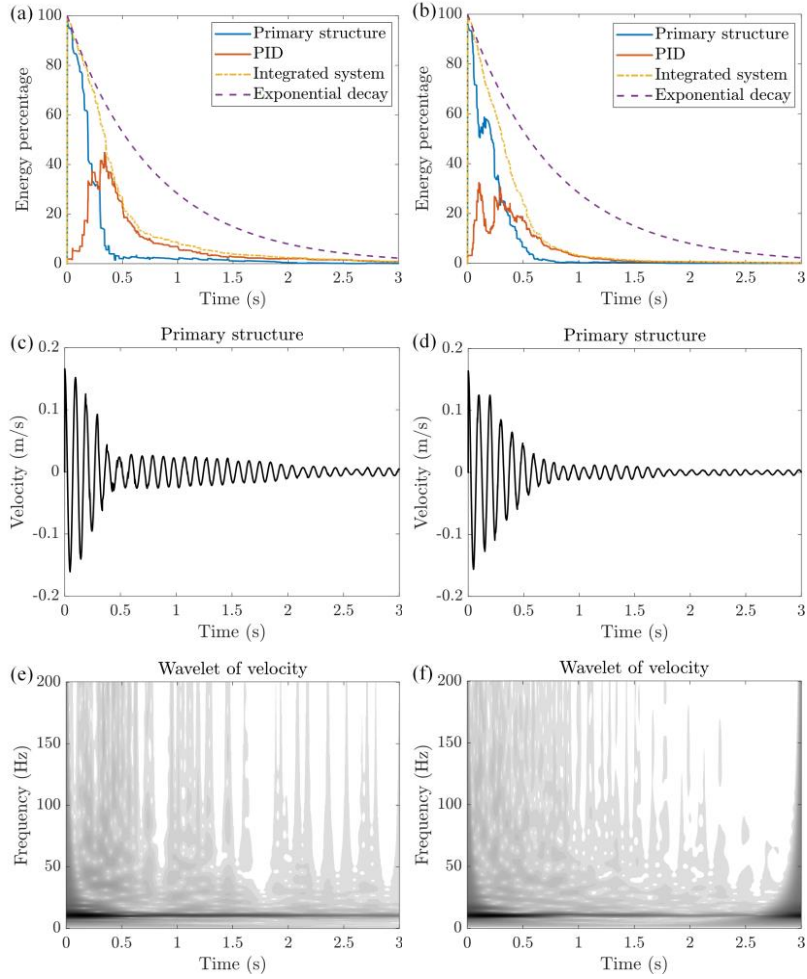

**Fig. 11.** Instantaneous normalized energy percentage (a,b), velocity of the PS (c,d) and its normalized wavelet transform spectrum (e,f) for $F_0 = 5 \times 10^3$ N: (a,c,e) PID with five granules, (b,d,f) PID with eight granules.

## 6. Conclusions

The nonlinear time-bandwidth characteristics and their tunability with energy of a shock-excited linear dissipative oscillator (denoted as the primary structure – PS) incorporating a particle impact damper (PID) were investigated. The PIDs possessed a varying number of granules (from one to eight) arranged in symmetric or asymmetric initial topologies, inside a rectangular rigid container. Following the application of the shock to the PS, due to initial clearances, granule-granule and granule-container inelastic collisions, as well as frictional effects due to relative granule rotations, yielded strongly nonlinear and highly discontinuous transient dynamics of the integrated PS-PID system. The inelastic granular collisions were modelled by a Hertzian contact law in the radial direction and Coulomb's tanh-law in the tangential direction. Moreover, the resulting highly complex granular dynamics were simulated by the discrete element method, while the convergence of the numerical results was carefully checked based on the criteria detailed in [3].

The nonlinear bandwidth, decay-time constant and time-bandwidth (T-B) product of the PS decaying response were computed based on the expressions of [35,36], as a way to assess the dissipative capacity of the considered PIDs. The main conclusions of this study are summarized as follows. For PIDs with the same number of granules, the initial granular topologies drastically affect the time-bandwidth results, and the time-bandwidth results are further highly tunable with energy. For the cases of optimal multiple-granule PIDs studied in [3], the nonlinear bandwidth of the PS is larger compared to the linear PID-less system over a broad shock intensity range, whereas the nonlinear bandwidth of the integrated PS-PID system is smaller than that of the linear system for small shock intensities and larger for large ones. However, based on the resulting T-B product, it was determined that for weak shocks the overall dissipative capacity of the integrated system is enhanced compared to the linear PID-less system with the same bandwidth as the integrated system, while it is diminished for stronger shocks. In general, the optimal PID configurations [3] are capable of effectively suppressing the shock-induced vibration of the PS over a broad energy range, while the integrated system can store shock energy for longer times at low shock levels.

High nonlinear bandwidth can be achieved by using multiple granules in optimal configurations [3], and this is evidenced by intense energy scattering from low-to-high frequencies together with strong irreversible energy transfer and rapid energy dissipation from the PS to the PID. Then the PID acts, in essence, as a granular nonlinear energy sink.

Due to inelastic granular collisions and frictional effects, the classical T-B limit can be violated in the considered systems. For optimal PIDs, the T-B product is below unity at low shock forces and becomes larger than unity at high shock levels; this is due to the fact that, as the shock intensity increases, the decay-time constant does not decrease in proportion to the increase of the nonlinear bandwidth. This directly affects the overall dissipative capacity of the integrated PS-PID system, which is enhanced compared to the linear PID-less system for weak shocks and diminished for strong shocks. This is one of the many results of this work that highlights the strong (passive) tunability of the performance of the PIDs with energy.

Assessing the dissipative capacity of a nonlinear resonator (even with discontinuous dynamics) based on its nonlinear bandwidth features provides a new way to evaluate the performance and capacity for shock mitigation of a rather broad class of engineering systems and structures. Given that the nonlinear bandwidth calculations are insensitive to the complexity and dimensionality of the dynamics, further applications are envisioned in areas such as seismic and blast mitigation, acoustic metamaterials, vibration energy harvesting, monitoring and sensing, and other fields.

**Acknowledgements**

This research was supported in part by the National Natural Science Foundation of China (No. 12202160), the China Scholarship Council (XL), and the Innovation Program of the Shanghai Municipal Education Commission (No. 2019-01-07-00-09-E00018). This support, which is gratefully acknowledged, made possible the academic visit of Xiang Li to the University of Illinois at Urbana-Champaign.

**Data availability**

The data generated and analyzed during the current study are available from the corresponding author on reasonable request.

**Declaration of competing interest**

The authors have no relevant financial or non-financial conflict of interest to disclose.

**References**

[1] Z. Lu, Z. Wang, S.F. Masri, X. Lu, Particle impact dampers: Past, present, and future, Struct. Control Heal. Monit. 25 (2018) 1–25. https://doi.org/10.1002/stc.2058.

[2] L. Gagnon, M. Morandini, G.L. Ghiringhelli, A review of particle damping modeling and testing, J. Sound Vib. 459 (2019) 114865. https://doi.org/10.1016/j.jsv.2019.114865.

[3] L. Xiang, A. Mojahed, L.Q. Chen, L.A. Bergman, A.F. Vakakis, Irreversible energy transfers in systems with particle impact dampers, Nonlinear Dyn. (accepted).

[4] M. Masmoudi, S. Job, M.S. Abbes, I. Tawfiq, M. Haddar, Experimental and numerical investigations of dissipation mechanisms in particle dampers, Granul. Matter. 18 (2016). https://doi.org/10.1007/s10035-016-0667-4.

[5] Z. Lu, B. Huang, Y. Zhou, Theoretical study and experimental validation on the energy dissipation mechanism of particle dampers, Struct. Control Heal. Monit. 25 (2018) 1–16. https://doi.org/10.1002/stc.2125.

[6] W. Xiao, Z. Chen, T. Pan, J. Li, Research on the impact of surface properties of particle on damping effect in gear transmission under high speed and heavy load, Mech. Syst. Signal Process. 98 (2018) 1116–1131. https://doi.org/10.1016/j.ymssp.2017.05.021.

[7] W. Yan, B. Wang, H. He, Research of mechanical model of particle damper with friction effect and its experimental verification, J. Sound Vib. 460 (2019) 114898. https://doi.org/10.1016/j.jsv.2019.114898.

[8] N. Meyer, R. Seifried, Damping prediction of particle dampers for structures under forced vibration using effective fields, Granul. Matter. 23 (2021). https://doi.org/10.1007/s10035-021-01128-z.

[9] Z. Lu, Y. Liao, Z. Huang, Stochastic response control of particle dampers under random seismic excitation, J. Sound Vib. 481 (2020) 115439. https://doi.org/10.1016/j.jsv.2020.115439.

[10] A. Sack, K. Windows-Yule, M. Heckel, D. Werner, T. Pöschel, Granular dampers in microgravity: sharp transition between modes of operation, Granul. Matter. 22 (2020) 1–6. https://doi.org/10.1007/s10035-020-01017-x.

[11] X. Ye, Y.Q. Ni, M. Sajjadi, Y.W. Wang, C.S. Lin, Physics-guided, data-refined

modeling of granular material-filled particle dampers by deep transfer learning, Mech. Syst. Signal Process. 180 (2022) 109437. https://doi.org/10.1016/j.ymssp.2022.109437.

[12] C. Zhang, Z. Zhao, T. Chen, H. Liu, K. Zhang, Discrete element method model of electromagnetic particle damper with a ferromagnetic end cover, J. Sound Vib. 446 (2019) 211–224. https://doi.org/10.1016/j.jsv.2019.01.034.

[13] N.R. Varela-Rosales, A. Santarossa, M. Engel, T. Pöschel, Granular binary mixtures improve energy dissipation efficiency of granular dampers, Granul. Matter. 25 (2023) 1–8. https://doi.org/10.1007/s10035-023-01337-8.

[14] C. Gnanasambandham, F. Fleissner, P. Eberhard, Enhancing the dissipative properties of particle dampers using rigid obstacle-grids, J. Sound Vib. 484 (2020). https://doi.org/10.1016/j.jsv.2020.115522.

[15] Y. Hu, H. Zan, Y. Guo, J. Jiang, Z. Xia, H. Wen, Z. Peng, Energy dissipation characteristics of particle dampers with obstacle grids, Mech. Syst. Signal Process. 193 (2023) 110231. https://doi.org/10.1016/j.ymssp.2023.110231.

[16] M. Żurawski, C. Graczykowski, R. Zalewski, The prototype, mathematical model, sensitivity analysis and preliminary control strategy for Adaptive Tuned Particle Impact Damper, J. Sound Vib. 564 (2023) 117799. https://doi.org/10.1016/j.jsv.2023.117799.

[17] X.M. Bai, L.M. Keer, Q.J. Wang, R.Q. Snurr, Investigation of particle damping mechanism via particle dynamics simulations, Granul. Matter. 11 (2009) 417–429. https://doi.org/10.1007/s10035-009-0150-6.

[18] A.F. Vakakis, O. V. Gendelman, L.A. Bergman, M.D. Michael, G. Kerschen, Y.S. Lee, Nonlinear Targeted Energy Transfer in Mechanical and Structural Systems, Springer Science & Business Media, Berlin, 2008.

[19] T. Yang, S. Hou, Z.H. Qin, Q. Ding, L.Q. Chen, A dynamic reconfigurable nonlinear energy sink, J. Sound Vib. 494 (2021). https://doi.org/10.1016/j.jsv.2020.115629.

[20] X.F. Geng, H. Ding, Theoretical and experimental study of an enhanced nonlinear energy sink, Nonlinear Dyn. 104 (2021) 3269–3291. https://doi.org/10.1007/s11071-021-06553-6.

[21] X.F. Geng, H. Ding, Two-modal resonance control with an encapsulated nonlinear energy sink, J. Sound Vib. 520 (2022) 116667.

https://doi.org/10.1016/j.jsv.2021.116667.

[22] H. Ding, Y. Shao, NES cell, Appl. Math. Mech. (English Ed. 43 (2022) 1793–1804. https://doi.org/10.1007/s10483-022-2934-6.

[23] Z.J. Wang, J. Zang, Y.W. Zhang, Method for Controlling Vibration and Harvesting Energy by Spacecraft: Theory and Experiment, AIAA J. 60 (2022) 6097–6115. https://doi.org/10.2514/1.J061998.

[24] Y. Zeng, H. Ding, A tristable nonlinear energy sink, Int. J. Mech. Sci. 238 (2023) 107839. https://doi.org/10.1016/j.ijmecsci.2022.107839.

[25] Y. Liu, Y. Wang, Vibration suppression of a linear oscillator by a chain of nonlinear vibration absorbers with geometrically nonlinear damping, Commun. Nonlinear Sci. Numer. Simul. 118 (2023) 107016. https://doi.org/10.1016/j.cnsns.2022.107016.

[26] J. Dou, H. Yao, Y. Cao, S. Han, R. Bai, Enhancement of bistable nonlinear energy sink based on particle damper, J. Sound Vib. 547 (2023) 117547. https://doi.org/10.1016/j.jsv.2022.117547.

[27] A.F. Vakakis, Passive nonlinear targeted energy transfer, Philos. Trans. R. Soc. A Math. Phys. Eng. Sci. 376 (2018) 20170132. https://doi.org/10.1098/rsta.2017.0132.

[28] B. Fang, T. Theurich, M. Krack, L.A. Bergman, A.F. Vakakis, Vibration suppression and modal energy transfers in a linear beam with attached vibro-impact nonlinear energy sinks, Commun. Nonlinear Sci. Numer. Simul. 91 (2020) 105415. https://doi.org/10.1016/j.cnsns.2020.105415.

[29] X. Li, A. Mojahed, L.Q. Chen, L.A. Bergman, A.F. Vakakis, Shock response mitigation of a large-scale structure by modal energy redistribution facilitated by a strongly nonlinear absorber, Acta Mech. Sin. 38 (2022) 121464.

[30] T. Theurich, A.F. Vakakis, M. Krack, Predictive design of impact absorbers for mitigating resonances of flexible structures using a semi-analytical approach, J. Sound Vib. 516 (2022). https://doi.org/10.1016/j.jsv.2021.116527.

[31] M. Gzal, B. Fang, A.F. Vakakis, L.A. Bergman, O. V. Gendelman, Rapid non-resonant intermodal targeted energy transfer (IMTET) caused by vibro-impact nonlinearity, Nonlinear Dyn. 101 (2020) 2087–2106. https://doi.org/10.1007/s11071-020-05909-8.

[32] J.R. Tempelman, A. Mojahed, M. Gzal, K.H. Matlack, O. V. Gendelman, L.A. Bergman,

A.F. Vakakis, Experimental Inter-Modal Targeted Energy Transfer in a cantilever beam undergoing Vibro-impacts, J. Sound Vib. 539 (2022) 117212. https://doi.org/10.1016/j.jsv.2022.117212.

[33] M. Gzal, J.E. Carrion, M.A. AL-Shudeifat, B.F. Spencer, J.P. Conte, A.F. Vakakis, L.A. Bergman, O. V. Gendelman, Seismic mitigation of a benchmark twenty-story steel structure based on intermodal targeted energy transfer (IMTET), Eng. Struct. 283 (2023) 115868. https://doi.org/10.1016/j.engstruct.2023.115868.

[34] A.F. Vakakis, O. V. Gendelman, L.A. Bergman, A. Mojahed, M. Gzal, Nonlinear targeted energy transfer: state of the art and new perspectives, Nonlinear Dyn. 108 (2022) 711–741. https://doi.org/10.1007/s11071-022-07216-w.

[35] A. Mojahed, L.A. Bergman, A.F. Vakakis, Generalization of the Concept of Bandwidth, J. Sound Vib. 533 (2022) 1–19. https://doi.org/10.1016/j.jsv.2022.117010.

[36] A. Mojahed, K.L. Tsakmakidis, L.A. Bergman, A.F. Vakakis, Exceeding the classical time-bandwidth product in nonlinear time-invariant systems, Nonlinear Dyn. 108 (2022) 3969–3984. https://doi.org/10.1007/s11071-022-07420-8.

[37] K.L. Johnson, Contact Mechanics, Cambridge University Press, Cambridge, 1985.

[38] Q. Zhang, W. Li, J. Lambros, L.A. Bergman, A.F. Vakakis, Pulse transmission and acoustic non-reciprocity in a granular channel with symmetry-breaking clearances, Granul. Matter. 22 (2020) 1–16. https://doi.org/10.1007/s10035-019-0982-7.

[39] C. Wang, Q. Zhang, A.F. Vakakis, Wave transmission in 2D nonlinear granular-solid interfaces, including rotational and frictional effects, Granul. Matter. 23 (2021) 1–21.

[40] C. Wang, S. Tawfick, A.F. Vakakis, Two-dimensional granular-thin plate interface for shock mitigation, Int. J. Non. Linear. Mech. 146 (2022) 104168. https://doi.org/10.1016/j.ijnonlinmec.2022.104168.

[41] N. Meyer, R. Seifried, Toward a design methodology for particle dampers by analyzing their energy dissipation, Comput. Part. Mech. 8 (2021) 681–699. https://doi.org/10.1007/s40571-020-00363-0.

[42] T.P. Sapsis, D.D. Quinn, A.F. Vakakis, L.A. Bergman, Effective stiffening and damping enhancement of structures with strongly nonlinear local attachments, J. Vib. Acoust. Trans. ASME. 134 (2012) 1–12. https://doi.org/10.1115/1.4005005.